\documentclass[12pt]{amsart}
\usepackage{amsfonts, amssymb,amsmath, amsthm, amsxtra, latexsym, amscd}
\usepackage[latin1]{inputenc} 

\def\draft{\centerline{(Draft {\the \day}/{\the\month} \the \year.)}}

\theoremstyle{definition}
\newtheorem{theo+}    {Theorem}      [section]
\newtheorem{prop+}  [theo+]  {Proposition}
\newtheorem{coro+}  [theo+]  {Corollary}
\newtheorem{lemm+}  [theo+]  {Lemma}
\newtheorem{deep+}  [theo+]  {Deep Result}
\newtheorem{fact+}  [theo+]  {Fact}
\theoremstyle{definition}
\newtheorem{exam+}  [theo+]  {Example}
\newtheorem{rema+}  [theo+]  {Remark}
\newtheorem{defi+}  [theo+]  {Definition}
\newtheorem{xca+}[theo+]{Exercise}

\numberwithin{equation}{section}

\def\draft{\centerline{(Draft {\the \day}/{\the\month} \the \year.)}}
\bigskip

\def\refn#1.#2{\expandafter\def\csname#1\endcsname{[#2]}}
\def\refnr#1.{\csname#1\endcsname}

\def\fa{\mathfrak a}

\def\ff{\mathfrak f}
\def\fg{\mathfrak g}
\def\fk{\mathfrak k}

\def\fp{\mathfrak p}

\def\a{\alpha}

\def\Claminv2{|C(\Lambda)|^{-2}}

\def\lam{\lambda}

\def\de{d\varepsilon}

\def\Aa2D{A^{\a,2}(D)}
\def\bAa2D{\overline{A^{\a,2}(D)}}
\def\Ab2D{A^{\beta,2}(D)}
\def\bAb2D{\overline{A^{\beta,2}(D)}}

\def\Norm#1_#2{\Vert#1\Vert_{#2}}

\def\phipl12{\phi_{p_{l_1}, p_{l_2}}}
\def\phip01{\phi_{p_{0}, p_{0}}}

\def\a{\alpha}

\def\Claminv2{|C(\Lambda)|^{-2}}

\def\sig{\sigma}
\def\Sig{\Sigma}

\def\lam{\lambda}

\def\ad{\operatorname{ad}}

\def\ch{\operatorname{ch}}
\def\coth{\operatorname{coth}}

\def\diag{\operatorname{diag}}
\def\min{\operatorname{min}}

\def\rank{\operatorname{rank}}
\def\tanh{\operatorname{tanh}}
\def\coth{\operatorname{coth}}
\def\ch{\operatorname{ch}}
\def\exp{\operatorname{exp}}

\def\rad{\operatorname{rad}}

\def\sh{\operatorname{sh}}
\def\tanh{\operatorname{tanh}}
\def\sh{\operatorname{sh}}

\def\tr{\operatorname{tr}}

\def\CH{\operatorname{CH}}
\def\SH{\operatorname{SH}}

\def\de{d\varepsilon}

\def\Aa2D{A^{\a,2}(D)}
\def\bAa2D{\overline{A^{\a,2}(D)}}
\def\Ab2D{A^{\beta,2}(D)}
\def\bAb2D{\overline{A^{\beta,2}(D)}}

\def\phipl12{\phi_{p_{l_1}, p_{l_2}}}
\def\phip01{\phi_{p_{0}, p_{0}}}

\def\bc{\mathbb C}
\def\br{\mathbb R}

\def\alg/{algebra} 
\def\Alg/{Algebra} 
\def\alt/{alternative} 
\def\anal/{analytic}
\def\analfunc/{\anal/\ \func/}
\def\Ans/{\it Answer. \normal}
\def\ass/{associative}
\def\nass/{non-\ass/}
\def\autom/{automorphism}
\def\homom/{homomorphism}
\def\isom/{isomorphism}
\def\bdd/{bounded}
\def\Bdd/{Bounded}
\def\bddsymdom/{bounded \sym/ \dom/}
\def\Cartdom/{Cartan \dom/}
\def\bdry/{boundary}
\def\bsd/{\bdd/ \symdom/}
\def\bv/{boundary value}
\def\cf/{{\it cf}\.}
\def\Cf/{{\it Cf}\.}
\def\charr/{character}
\def\coeff/{coefficient}
\def\comm/{commutative}
\def\cpct/{compact}
\def\compl/{complex}
\def\comp/{complex}
\def\Comp/{Complex}
\def\conf/{conformal}
\def\conj/{conjugate}
\def\conn/{connect}
\def\cont/{continuous}
\def\conv/{converge} 
\def\convc/{convergence}
\def\convt/{convergent}
\def\convx/{convex}
\def\coord/{coordinate}
\def\lcoord/{local coordinate}
\def\Corr/{Corresponding}
\def\corr/{corresponding}
\def\corrd/{correspond}
\def\cov/{covariant}
\def\decomp/{decomposition}
\def\deco/{decompose}
\def\diff/{different} 
\def\Diff/{Different} 
\def\dimn/{dimension} 
\def\distr/{distribution} 
\def\div/{diverge} 
\def\dom/{domain}
\def\eg/{\hbox{\it e.g}\.}
\def\eigenf/{eigen\-\func/}
\def\eigensp/{eigen\-space}
\def\eigenv/{eigen\-value}
\def\eq/{equation}
\def\equa/{equation}
\def\de/{\diff/ial \equa/}
\def\do/{\diff/ial operator}
\def\ode/{ordinary \de/}
\def\pde/{partial \de/}
\def\pdo/{partial \diff/ial operator}
\def\psdo/{pseudo \diff/ial operator}
\def\fin/{finite}
\def\Ex/{\it Example.\ \normal}
\def\Exnr#1/{\it Example #1.\ \normal}
\def\foll/{follow}
\def\follg/{following}
\def\Follg/{Following}
\def\func/{function}
\def\Func/{Function}
\def\Fonc/{Fonc\-tion}
\def\fonc/{fonc\-tion}
\def\Funk/{Funk\-tion}
\def\funk/{Funk\-tion}
\def\gen/{general}
\def\har/{harmonic}
\def\Hint/{\it Hint. \normal}
\def\hist/{historic}
\def\histcl/{historical}
\def\hol/{holo\-morphic}
\def\homog/{ho\-mo\-ge\-ne\-ous}
\def\hyp/{hyper\-bolic}
\def\hyperg/{hyper\-geometric}
\def\ie/{\hbox{\it i.e.}}
\def\iff/{if and only if}
\def\ineq/{inequality}
\def\infra/{{\it inf\-ra}}
\def\ultra/{{\it ult\-ra}}
\def\Inpart/{In particular}
\def\inpart/{in particular}
\def\instof/{instead of}
\def\interps/{interpolation space}
\def\interp/{interpolation}
\def\Interp/{Interpolation}
\def\interpr/{Interpretation}
\def\Intr/{Introduction}
\def\intv/{interval}
\def\inv/{invariant}
\def\invc/{invariance}
\def\Iowords/{In other words}
\def\iowords/{in other words}
\def\ipr/{inner product}
\def\irred/{irreducible}
\def\lb/{line bundle}
\def\lin/{linear}
\def\lhs/{left hand side}
\def\rhs/{right hand side}
\def\loc/{local}
\def\math/{mathematic} 
\def\mathcn/{\math/ian}
\def\manif/{manifold}
\def\meas/{measure}
\def\measl/{measurable}
\def\mero/{mero\-morphic}
\def\mon/{monomial}
\def\monog/{monogenic}
\def\mult/{multiple}
\def\multy/{multiply}
\def\multn/{multiplication}
\def\nas/{necessary and sufficient}
\def\nbd/{neighborhood}
\def\neg/{negative}
\def\nondeg/{nondegenerate}
\def\Oohand/{On the other hand}
\def\oohand/{on the other hand}
\def\Oonhand/{On the one hand}
\def\oonhand/{on the one hand}
\def\oper/{operator}
\def\orth/{ortho\-gonal}
\def\orthon/{ortho\-normal}
\def\otoh/{on the other hand}
\def\quat/{quaternion}
\def\pp/{\hbox{a. e.}}
\def\psh/{plurisubharmonic}
\def\pol/{polynomial}
\def\pot/{potential}
\def\pos/{positive}
\def\princ/{principle}
\def\prob/{probability}
\def\proj/{projective}
\def\projn/{projection}
\def\Proof/{\it Proof:\normal}
\def\Rem/{\it Remark\normal}
\def\Remnr#1/{\it Remark\ \normal #1. }
\def\rep/{representation}
\def\meta/{metaplectic representation}
\def\repr/{reproducing}
\def\reprker/{reproducing kernel}
\def\resp/{respective} 
\def\resply/{respectively}
\def\restr/{restriction}
\def\sa/{self-adjoint}
\def\st/{such that}

\def\sol/{solution}
\def\ru/{space}
\def\sph/{spherical}
\def\ssp/{sub\ru/}
\def\sym/{symmetric}
\def\Sym/{Symmetric}
\def\symb/{symbol}
\def\symbc/{symbolic}
\def\symdom/{\sym/ domain}
\def\symp/{symplectic}
\def\Theor#1/{\fet Theorem #1.\ \normal}
\def\Lem#1/{\fet Lemma #1.\ \normal}
\def\Lemma/{\fet Lemma.\ \normal}
\def\topl/{topology}
\def\topll/{topological}
\def\transf/{transform}
\def\transl/{translation}
\def\transfn/{transformation}
\def\transv/{transvectant}
\def\trig/{trigonometric}
\def\tril/{trilinear}
\def\trilf/{trilinear form}
\def\uhp/{upper halfplane}
\def\uhs/{upper halfspace}
\def\vb/{vector bundle}
\def\vf/{vector field}
\def\vsp/{vector space}
\def\wrt/{with respect to}
\def\Wlog/{Without loss of generality}
\def\a{\alpha}

\def\lam{\lambda}

\def\sig{\sigma}
\def\Sig{\Sigma}

\def\Ab/{Abel}
\def\Ban/{Banach}
\def\Bansp/{\Ban/ space}
\def\Belt/{Bel\-tra\-mi}
\def\Berg/{Berg\-man}
\def\Bern/{Ber\-nou\-lli}
\def\Berz/{Berezin}
\def\Bess/{Bessel}
\def\Cart/{Car\-tan}
\def\Cay/{Cay\-ley}
\def\CG/{Clebsch-Gordan}
\def\Cl/{Clifford}
\def\CR/{Cauchy-Rie\-mann}
\def\Dir/{Dirichlet}
\def\Eucl/{Euclide}
\def\F/{Fourier}
\def\Hank/{Hankel}
\def\Hankf/{\Hank/ form}
\def\Herm/{Hermite}
\def\Hilb/{Hilbert}
\def\Hilbs/{Hilbert space}
\def\Hilbsp/{Hilbert space}
\def\HS/{Hilbert-Schmidt}
\def\Lag/{La\-grange}
\def\Lap/{La\-place}
\def\LapBelt/{\Lap/-\Belt/}
\def\Leb/{Lebesgue}
\def\Marc/{Mar\-cin\-kie\-wicz}
\def\Moeb/{Moebius}
\def\Moebt/{Moebius transformation}
\def\Moebtransfn/{Moebius transformation}
\def\Pla/{Plan\-che\-rel}
\def\Poin/{Poin\-car\'e}
\def\Riem/{Rie\-mann}
\def\Riemn/{\Riem/ian}
\def\psRiemn/{pseudo-\Riem/ian}
\def\Riems/{Rie\-mann surface}
\def\Schroe/{Schr\"odinger}
\def\Weier/{Weier\-strass}
\def\anal/{analytic}
\def\bsd/{bounded symmetric domain  }
\def\bdd/{bounded}
\def\calc/{calculation}\def\conj{conjugate}
\def\calci/{calculating}\def\eg{e.g.}
\def\conj/{conjugate}
\def\deco/{decomposition}
\def\eg/{e.g.}
\def\fct/{function}
\def\gp/{group}
\def\hw/{highest weight}
\def\hwv/{highest weight vector}
\def\hwvs/{highest weight vectors}
\def\lw/{lowest weight}
\def\lwv/{lowest weight vector}
\def\lwvs/{lowest weight vectors}
\def\hds/{holomorphic discrete series}
\def\iff/{if and only if}
\def\inv/{invariant}
\def\irrde/{irreducible decomposition}
\def\meas/{measure}
\def\transf/{transform}
\def\rep/{representation}
\def\resp/{respectively}
\def\inters/{intertwines}
\def\interg/{intertwining}
\def\meta/{metaplectic representation}
\def\qu/{quaternion}
\def\rep/{representation}
\def\symdom/{ symmetric domain}
\def\st/{such that}
\def\shd/{subhead}
\def\transf/{transform}
\def\wrt/{with respect to}

\def\Norm#1#2#3{\Vert#1\Vert^{#3}_{{#2}¥}}

\def\tr{\operatorname{tr}}

\begin{document}
\title[Radon transform]
{Radon transform on  symmetric matrix domains}
\author{Genkai Zhang}
\address{Department of Mathematics, Chalmers University of Technology and
G\"oteborg University
 , G\"oteborg, Sweden}
\email{genkai@math.chalmers.se}

\thanks{Research  supported by the Swedish
Science Council (VR)}
\subjclass{22E45, 33C67, 43A85, 44A12}
\keywords{Radon transform,   inverse Radon transform, symmetric domains,
Grassmannians  manifolds,  Lie groups,
fractional integrations, Bernstein-Sato formula, Cherednik operators,
invariant differential operators}

\def\bbK{\mathbb K}
\def\bKk{\mathbb K^k}
\def\bKn{\mathbb K^n}
\def\Cos^2{\text{Cos}^2}

\begin{abstract}
Let $\bbK=\mathbb R, \mathbb C, \mathbb H$ 
be   the  field of real, complex or quaternionic numbers
and $M_{p, q}(\bbK)$ the vector space of all
$p\times q$-matrices.
Let $X$ be
 the matrix unit ball in
$M_{n-r, r}(\bbK)$ consisting of contractive
matrices. As a symmetric space,
 $X=G/K=O(n-r, r)/O(n-r)\times O(r)$,
$U(n-r, r)/U(n-r)\times U(r)$ and
respectively  $Sp(n-r, r)/Sp(n-r)\times Sp(r)$.
The matrix unit ball $y_0$ in $M_{r^\prime-r, r}$ with $r^\prime \le n-1$
is a totally geodesic submanifold
of $X$ and let $Y$ be the set of  all
 $G$-translations of the submanifold $y_0$. The set $Y$
is then a manifold and an affine symmetric space.
We consider the Radon transform $\mathcal Rf(y)$ for
 functions $f\in C_0^\infty(X)$ defined
by integration of $f$ over the subset $y$,
and the dual transform $\mathcal R^t F(x), x\in X$ 
for functions $F(y)$ on $Y$.
For $2r <n, 2r\le r^\prime$ with certain evenness condition
in the case $\bbK=\br$,  we find an $G$-invariant 
differential operator $\mathcal M$ 
and prove it is the right inverse
of $\mathcal R^t \mathcal R$, 
$\mathcal R^t \mathcal R \mathcal M f=c f$, for
$f\in C_0^\infty(X)$, $c\ne 0$.
The operator $f\to \mathcal R^t\mathcal Rf$ 
is an integration of $f$ against
a (singular)  function determined by
the root systems of $X$ and $y_0$.
We study the analytic continuation
of the powers of the function and we find
a Bernstein-Sato type formula
generalizing earlier work of the author
in the set up of Berezin transform. When $X$
is rank one domain of hyperbolic balls in 
$\mathbb K^{n-1}$ and $y_0$
is the hyperbolic ball in $\mathbb K^{r^\prime -1}$,
$1<r^\prime<n$ we obtain an inversion
formula for the Radon transform, namely $\mathcal M\mathcal R^t\mathcal R f=c f$.
This  generalizes  earlier results
of Helgason for non-compact rank one symmetric
spaces for the case $r^\prime=n-1$.

\end{abstract}

\maketitle








\def\mP{\mathcal P}
\def\mF{\mathcal F}
\def\SdN{\mathcal S^N}
\def\eSdN{\mathcal S_N^\prime}
\def\Gc{G_{\mathbb C}}
\def\rqs{G\backslash K_{\mathbb C}}

\baselineskip 1.25pc

\section{Introduction}

Let $\bbK$ be the field
of the real, complex and quanternionic numbers.
In this paper we will study the Radon transform on non-compact
symmetric
domains of $(n-r)\times r$-matrices over $\bbK$.
To motivate and formulate our results recall
first the Radon transform on compact Grassmannian manifolds.
Consider the Grassmannian manifold
$X_c=G(n, r)$ of $r$-dimensional subspaces
in $\bbK^n$. Let $1\le r <r^\prime < n$. The Radon transform 
$\phi(y)=\mathcal Rf(y)$, $y\in G(n, r^\prime)$, of  $f(x)$ on $G(n, r)$
is defined as  the  integration over the subset of
all subspaces
$x\in G(n, r)$ which are contained in the fixed subspace
$y\in Y_c=G(n, r^\prime)$. The dual Radon transform $(\mathcal R^t \phi)(x)$
for a   function $\phi(y)$ on $G(n, r^\prime)$ is defined
as an integration over the subset of all $y\in G(n, r^\prime)$ which 
contain the given subspace $x\in G(n, r)$.
One natural
question is to invert the Radon transform. 
To do that one observes first that
the space $X_c=G(n, r)$ is a compact symmetric space
$X_c=G_c/K$, $G_c=U(n; \bbK)$, $K=U(r; \bbK)\times U(n-r; \bbK)$
and that the operator $\mathcal R^t \mathcal R$ defined on  $X_c$
and is $G_c$-invariant. The $L^2$-space $L^2(X_c)$ of functions on $X_c$
is decomposed as a direct sum of irreducible subspaces
of $G_c$ with multiplicity free by the Cartan-Helgason theorem.
In his paper \cite{Grinberg-grs}
Grinberg computes the eigenvalues of $\mathcal R^t \mathcal R$
on the irreducible subspaces. (More precisely the operator
$\mathcal R$ maps highest weight vectors to highest weight vectors
in the respective functions spaces on $X_c$ and on $Y_c$,
and Grinberg computes the eigenvalue of $\mathcal R$
under some suitable normalization, which then gives
 the eigenvalues of $\mathcal R^t \mathcal R$.)
For the
complex,  quanternionic  and real Grassmannian
with $\rank(X)\le \rank(G(n, r^\prime))$
and  certain even conditions in the real case
it turns out
that there exists an invariant differential
operator $\mathcal M$ which is the left inverse of
 $\mathcal R^t \mathcal R$, namely
$\mathcal M \mathcal R^t \mathcal Rf =f$ for smooth functions  $f$;
thus $\mathcal M \mathcal R^t $ is the left inverse of $\mathcal R$
 giving an inversion formula for $\mathcal R$.
Some different inversion formulas
had also been studied earlier,
see \cite{Gelfand-Graev-Rosu}
 and references therein.

Recently
some explicit constructions of the invariant
differential operator have been studied, see e. g. 
\cite{Kakehi-jaf99}
and \cite{Gonzales-Kakehi-1}.
It is a natural and very interesting question
to invert the Radon transform on non-compact
symmetric spaces, namely replacing
the compact spaces $X_c$ and $Y_c$
by their non-compact duals
$X$ and $Y$.

When $r=1$ the space
$G(n, 1)$ is  the projective space
 $P^{n-1}(\bbK)$, and it non-compact dual
space   $X=G/K$ can be realized as the unit ball in $\bbK^{n-1}$.
The unit ball $y_0$ in $\bbK^{r^\prime-1}$ can  be
realized as a totally geodesic submanifold of $X$.
We let $Y$ be the set of all $G$-translates
of $y_0$; $Y$ is then an affine symmetric space
and is in particular an analytic manifold.
In  \cite[Chapter I]{He1} Helgason 
 studied
the Radon transform from functions on $X$ to $Y$
for the case when $r^\prime =n-1$ 
and has found an inversion formula by a case
by case computation;
  in the special case of real hyperbolic spaces
in $\br^{n-1}$
he has also found
an inversion formula for any $r^\prime \le n-1$,
see \cite[Chapter III, Theorem 1.5]{He-radonbk},
and \cite{Rubin-adv-math}.

 Radon transform can also be defined
on higher rank non-compact symmetric spaces
by integration over totally geodesic submanifolds, but
finding its inversion has become rather complicated, by several reasons. 
Firstly the  decomposition of $L^2(X)$ under $G$
 is a direct integral (rather than the direct sum)   and 
 the  operator
 $\mathcal R^t\mathcal R$ is $G$-invariant but  not  bounded,
so that integration against the spherical functions is analytically
complicated, compared with the compact
case \cite{Grinberg-grs}. Secondly,
 in contrast to the rank one case
\cite[Chapter I]{He1}
 it is  not clear
if one can  find
an invariant differential operator as a left inverse of the operator $\mathcal R^t\mathcal R$,
because one does not  know if operator $\mathcal R^t\mathcal R$ commutes
with invariant differential operators (and we believe it is
not); if, on the other hand, we want
to find a right inverse of  $\mathcal R^t\mathcal R$,
then the  analytic continuation of
the integration defining $\mathcal R^t\mathcal R$
is difficult to handle since generally speaking no explicit
construction of invariant differential operators
is given.

In the present paper we will study Radon transform
on the non-compact dual $X=G/K$ of the Grassmannian
manifold $X_c=G_c/K$, with $G=U(n-r, r; \bbK)$. The space $X$ can be realized
as the matrix unit  ball of $(n-r)\times r$-matrices.
There is a natural totally geodesic submanifold $y_0$
consisting of  $(r^\prime-r)\times r$-matrices
in $X$. The set of all $G$-translates
of $y_0$
is a manifold $Y$, and is an affine symmetric space
$$
Y=G/H, \quad H=U(n-r^\prime; \bbK)\times U(r^\prime-r, r; \bbK).
$$
of rank $\min\{r^\prime, n-r^\prime\}$.
We assume 
\begin{equation}
  \label{ran-cd}
2r\le \min\{r^\prime, 2(n-r^\prime)\},
\end{equation}
so that the submanifold $y_0$ is
of the same rank as $X$, which in turn
is less than that of $Y$.

The transform
$\mathcal R^t \mathcal R$ on the space $C_0^\infty(X)$
is  $G$-invariant. 
 The value of $\mathcal R^t \mathcal Rf(x)$
at $x=o$ is an integration of $f$
against certain power of the
hyperbolic sine function $|\SH(t)|^\delta$
(see the formula  ( \ref{eq:sinht}))
determined by the root systems
of $y_0$ and $X$.
We study thus the meromorphic continuation
 $|\SH(t)|^\delta$ in $\delta\in \mathbb C$. For that
purpose we find certain 
Bernstein-Sato type formula (see Theorem \ref{BSF})
for the  function $|\SH(t)|^\delta$ using the Cherednik
operators for general root multiplicities, which
generalizes an earlier formula
for the function $|\CH(t)|^\delta$ in the context of
Berezin transform (\cite{gz-imrn}, \cite{gkz-manu-mat})
and which might be of independent interest. 
 The meromorphic continuation
turns out to be rather complicated, 
and by using our Bernstein-Sato
formula and earlier work of Gindikin
on Riesz type distributions
we are able to prove that  by proper normalization
it has analytic continuation and that
its gives the Dirac delta-distribution
at the base point; 
 see Theorem 3.4 and Proposition \ref{prop-dirac}.
We can then find the right inverse
of the transform $\mathcal R^t\mathcal R$, proving
that
$$
\mathcal R^t\mathcal R\mathcal M f= cf
$$
where $\mathcal M$ is an invariant
differential operator on $X$ given 
in terms of the Cherednik operators
(and we can thus find its eigenvalues
on the spherical functions, namely
its image under the Harish-Chandra homomorphism).

In  the case of rank one domains mentioned above
with $X$ being the unit  ball in $\mathbb K^{n-1}$
and  $y_0\subset X$ the ball in  $\mathbb K^{r^\prime-1}$.
We can
also include the exceptional rank one
domain into our consideration, namely $X$
being the unit ball in $\mathbb O^2$
and $y_0$  the unit ball in $\mathbb O=\mathbb R^8$.
In this case we can find an inversion formula
for the Radon transform due to
the commutativity of $\mathcal M$ with $\mathcal R^t \mathcal R$.
We give a systematic  treatment for all cases.
So our result
generalizes that of Helgason for $r^\prime =n-1$.  
To our knowledge our results are  new 
also in the rank one case.

I would like to thank the organizers of the workshop ``Complex Analysis, Operator Theory
and Mathematical Physics'' at Erwin Schr\"ordinger institute
for the invitation and  I thank
the institute for its support and hospitality.
I thank also Professors C. Dunkl, H. Schlichtkrull and H. Upmeier
for some helpful correspondences and discussions.

\section{Symmetric Matrix Domains and Radon transform}

In this section we introduce
the symmetric domain $X$, the Radon transform
$\mathcal R$ and the dual Radon transform $\mathcal R^t$.

\subsection{The matrix domain $X$}

The matrix domain $X$ can be introduced in the
general framework of Jordan triples \cite{Loos-bsd},
or symmetric spaces (\cite{He1} and \cite{He2}). We recall
briefly some necessary and elementary results. 

Let $\mathbb K=\br, \bc, \mathbb H$ be  the
field of real, complex or quaternionic numbers
with real dimension 
$$
a:=\dim_{\mathbb R} \bbK=1, 2, 4
$$
and let $x\to \bar x$ be the standard conjugation.
Let $M_{l, k}:=M_{l,k}(\mathbb K)$ be the  space
of  $n\times k$-matrices with entries $x$ in $\mathbb K$.
For brevity we denote (with some abuse
of notation) by $U(l, k)=U(l, k; \bbK)$
 the group of $\mathbb K$-linear transformations $g$
of $\mathbb K^{l+k}$ preserving the following
quadratic 
form  in $\mathbb K^{l+k}$,
$$\bar x_1 x_1+\cdots+ \bar x_l x_l -
\bar x_{l+1} x_{l+1}-\cdots- \bar x_{l+k} x_{l+k}.
$$
Here by a $\bbK$-linear transformation $g$ we mean
$g(vc)=g(v)c$, $v\in \bbK^k$, $c\in\bbK$.
Then $U(l, k)$ is  the classical group
$O(l, k), U(l, k), Sp(l, k)$ accordingly.

We consider, for $1\le r <n$, 
the domain
\begin{equation}
  \label{eq:def-X}
X=\{x\in M_{n-r, r}; I_r - x^\ast x >0\}.
\end{equation}
Here by $I_r - x^\ast x >0$ we mean
that $I_r - x^\ast x $ is a positive
$\mathbb R$-linear transform of the
real Hilbert space $\bbK^r=\mathbb R^{dr}$
(with the Euclidean metric 
$x\to x^\ast x$).

We fix the zero $(n-r)\times r$ matrix 
as a  reference point in $X$ and denote it by $o$.
The domain $X$  is a realization of the symmetric
space $G/K$,
$$
X=G\cdot o= G/K, \quad G=U(n-r,r), \quad K=U(n-r)\times U(r).
$$
The group action will be denoted by $gx$ or $g\cdot x$.
The domain $X$ can also be
characterized using the Jordan triple structure
$\{x\bar yz\}=xy^\ast z +z y^\ast x$
of $M_{n-r, r}$;  see \cite{Loos-bsd}.

Let $\fg$ be the Lie algebra of $U(n-r, r)$
and $\fg=\fk+\fp$  the Cartan decomposition. The space 
$\fp$ consists of $n\times n =(n-r +r)\times (n-r +r)$
block matrices of the form
$$
\begin{bmatrix} 0_{n-r, n-r}&x \\
                x^\ast &0_{r, r}\end{bmatrix}, \quad x\in M_{n-r, r},
$$
and each of them   will be identified with the upper left block
matrix $x\in M_{n-r, r}$. Here
$0_{k, l}$ stands for the zero $k\times l$-matrix.

Assume  $r\le n-r$.  The space $X$ is then of rank $r$.
Let
$\fa\subset \fp=M_{n-r, r}$ be the subspace consisting
of matrices 
$x$
of the form
$$
x=\begin{bmatrix} 0_{n -2r, r}\\
\diag\{t_1, t_2, \cdots,  t_r\}
\end{bmatrix}
=t_1\xi_1 +\cdots +t_r\xi_r
, \quad t_j\in \br 
$$
where $\diag\{t_1, t_2, \cdots,  t_r\}$ is the
diagonal $r\times r$-matrix with diagonal
entries $t_1, t_2, \cdots, t_r$ and
where the basis vectors $\xi_1, \cdots, \xi_r$ are self-defined
by the formula. We identify $\fa$ with $\br^r$
under the basis. Let $\{\xi_j^\ast\}$ be the 
basis of $\fa^\ast$ dual to $\{\xi_j\}$,
$\xi_j^\ast(\xi_j)=\delta_{j, k}$.
The root system $\Sig(\fg, \fa)$ of  $\fa$ in $\fg$ 
is of type $B$ or $D$ if $\bbK=\br$
or type $BC$ if $\bbK=\bc, \mathbb H$,
\begin{equation}
\label{rt-sys}
\Sig(\fg, \fa)=\{{\pm\xi_j^\ast \pm \xi_k^\ast}, j\ne k\} 
\cup\{\pm \xi_k^\ast\}
\cup\{\pm 2\xi_k^\ast\}
\end{equation}
with multiplicities $a$, $a(n-2k)$ and $a-1$ for
the respective sets of roots.
It is understood that the roots with
multiplicity zero will not appear
(e.g if $a=1$  the third set is empty).

We denote by $W$ the Weyl group of the root system
and fix an ordering
so that  $\xi_1^\ast >\cdots >\xi_r^\ast$
for type $B$  or $D$ and 
$\xi_1^\ast >\cdots >\xi_r^\ast >0$
for type BC.

Let $A=\exp(\fa)\subset G$. Then $A\cdot o$ is
a flat submanifold of $X$,
$$
A\cdot o=\{\begin{bmatrix}
0_{n-2r, r}\\ \diag\{\tanh t_1, \cdots, \tanh t_r\}
\end{bmatrix}; t\in \br^r\}
$$
Under the Cartan decomposition $X=KA\cdot o$  the
$G$-invariant Riemannian measure on $X$ is
(up to  constants) 
given by
\begin{equation}
  \label{cartan-x}
\int_X f(x) d\mu(x)=
\int_K \int_{\mathbb R^r}f(k\exp(\exp(t_1\xi_1+\cdots +t_r \xi_r)\cdot o
) )
d\mu(t) dk  
\end{equation}
where (some slight abuse of the notation $d\mu$)
\begin{equation}
\begin{split}
\label{rdl-ms}
d\mu(t)&=\prod_{\a\in \Sig_+}|2\sh (\a(t))|^{k_{\a}}dt\\
&=\prod_{j=1}^r |2\sh (2t_j)|^{a-1}
|2\sh (t_j)|^{a(n -2r)} 
 \prod_{i<j} |2\sh(t_i)\pm \sh(t_j)|^{a}dt_1\cdots dt_r,
\end{split}
\end{equation}
and $dk$ is the normalized Haar measure on $K$;
see \cite[Chapter I, Section 5]{He2}.

\subsection{The space $Y$ of totally geodesic submanifolds}
We assume throughout the rest of the paper that
the rank condition (\ref{ran-cd}) holds.
We will view the matrix space 
$M_{r^\prime -r, r}$ as a subspace of 
$M_{n-r, r}$ by identifying each $x\in 
M_{r^\prime -r, r}$ with the matrix
$\begin{bmatrix} 0\\x\end{bmatrix}$
in $M_{n-r, r}$.

Let $y_0$ be the following submanifold
of $X$,
\begin{equation}
  \label{eq:def-y-0}
y_0=X
\cap M_{r^\prime -r, r}=\{\begin{bmatrix} 0\\x\end{bmatrix}; 
x\in M_{r^\prime -r, r}, x^\ast x <I_{r}\}.  
\end{equation}
Then $y_0$ is a totally geodesic submanifold
of $X$ since $M_{r^\prime -r, r}$ is 
a Jordan subtriple  of $M_{n -r, r}$, see
\cite{Loos-bsd}. It is itself
a symmetric domain
\begin{equation}
  \label{eq:ss-y-0}
y_0=G_0/K_0, \quad G_0=U(r^\prime -r, r), \quad 
K_0=U(r^\prime -r)\times U(r)
\end{equation}
where $U(r^\prime -r, r)$ is also realized as
a subgroup of $G$ by identifying $g\in 
U(r^\prime -r, r)$ with $\diag \{I_{n-r^\prime}, g\}\in G$. 
By our assumption (\ref{ran-cd}) $y_0$ is also of rank $r$. Let
$\fg_0=\fk_0+\fp_0$ be the Cartan decomposition
of the Lie algebra $\fg_0$ of $G_0$. The space
$\fa\subset \fp$ is also a maximal abelian subspace of $\fp_0$
and the root system $\Sig(\fg_0, \fa)$ is
of the same type as that of $\Sig(\fg, \fa)$ in (\ref{rt-sys}) 
as a union of three sets,
\begin{equation}
\label{rt-sys-g0}
\Sig(\fg_0, \fa)=\{{\pm\xi_j^\ast \pm \xi_k^\ast}, j\ne k\} 
\cup\{\pm \xi_k^\ast\}
\cup\{\pm 2\xi_k^\ast\}
\end{equation}
 with respective
multiplicities $a, a(r^\prime-2r), a-1$.

We let 
\begin{equation}
  \label{eq:def-Y}
Y=G\cdot y_0
=\{gy_0\subset D; g\in G\}  
\end{equation}
be the set of all $G$-translations of $y_0$.
It is elementary to check that the stabilizer
of $y_0$ is 
$$
H=U(n-r^\prime)\times G_0 = U(n-r^\prime)\times  U(r^\prime-r, r).
$$
Thus $Y=G/H$
 is an affine symmetric space of rank $\min\{r^\prime, n-r^\prime\}$
and is in particular an analytic manifold.

\subsection{Radon transform}

We first  fix the normalization of the $G_0$-invariant
measure $d\mu_0$ on the Riemannian symmetric space $y_0$.
Under the  Cartan decomposition
of $y_0$, $y_0=K_0 A \cdot o$,
we have 
\begin{equation}
  \label{eq:nt-yo}
\int_{y_0}f(x)d\mu_0(x)
=\int_{K_0}\int_{\fa} f(k\exp(t_1\xi_1 +\cdots +t_r\xi_r)\cdot o) d\mu_0(t)dk  
\end{equation}
where $d\mu_0(t)$ is given  as in
(\ref{rdl-ms}) 
for $d\mu(t)$ with the multiplicity $a(n -2r)$ is replaced by
$a(r^\prime -2r)$, viz.
\begin{equation}
\label{rdl-ms-y0}
d\mu_0(t)=\prod_{j=1}^r |2\sh (2t_j)|^{a-1}
|2\sh (t_j)|^{a(r^\prime -2r)} 
 \prod_{i<j} |2\sh(t_i)\pm \sh(t_j)|^{a}dt_1\cdots dt_r.
\end{equation}
For any $y=g\cdot y_0\in Y$ there is a corresponding Riemannian
measure $d\mu_y$
on $y$ via the isometric mapping $g: y_0 \to y$; 
it is independent of the representatives $g$
by the  invariance by $d\mu_0$ under the stabilizer $H$ of $y_0$.

We define the Radon transform of $\mathcal R f$
for a function $f$ on $X$, by
\begin{equation}
\label{Rad-tr}
\mathcal Rf(y)=\int_{y}f(x)d\mu_y(x), \quad y\in Y
\end{equation}
whenever the integral is absolute convergent. 

To define the dual 
Radon transform we claim first that
$$
\{y\in Y; o\in y\}=K\cdot y_0:=\{k\cdot y_0; k\in K\}.
$$
 Indeed, it is clearly
that  $o=k\cdot o\in  k\cdot y_0$
 for any
$k\in K$ so that $K  \cdot y_0\subset
\{y\in Y; o\in y\}$. Conversely suppose $o\in y$. Write $y=g\cdot y_0$.
Thus  $o=g z$ for some  $z\in y_0$ which
in turn can be written as $z=g_z \cdot o$
for $g_z\in G_0$. We have then  $o=gz=g g_z \cdot o$,
and  that $g g_z=k$ for some $k\in K$, namely
$g=kg_z^{-1}$ and $y=g\cdot y_0=kg_z^{-1}\cdot y_0
=k y_0\in K\cdot y_0 $, proving the claim.
By the homogeneity of $X=G/K$ we
have
\begin{lemm+} For each $x\in X$ 
we have
$$
Y_x:=\{y\in Y; x\in y\}=g_x K \cdot y_0:=\{g_xk\cdot y_0; k\in K\},
$$
where $g_x\in G$ is such
that $x=g_x\cdot o$.
\end{lemm+}
The above lemma can also be explained by using the
double filtration in \cite{He2}
In terms of the above lemma
we can define the dual Radon transform for a smooth function 
$\psi$ on $Y$ by
\begin{equation}
  \label{eq:dual-radon}
\mathcal R^t \psi (x)=\int_{K}\psi(g_xky_0)dk, \quad x\in X  
\end{equation}
where $x=g_x\cdot o$. The function $\mathcal R^t \psi$ is clearly
well-defined since it's an integration a smooth function
over  compact manifolds,
and it can also be
written as 
\begin{equation}
  \label{eq:dual-radon-var}
\mathcal R^t \psi (x)=\int_{Y_x}\psi(y)d\nu_x(y), \quad x\in X  
\end{equation}
where $d\nu_x(y)$ is some probability measure on $Y_x$.

The operator $\mathcal R^t\mathcal R: C_0^\infty(X)\to
C^\infty(X)$ is then
\begin{equation}
  \label{eq:square-radon}
(\mathcal R^t \mathcal R) f (x)=\int_{K}\int_{y_0}
f(g_x k \eta)d\mu_0(\eta)dk, \quad x=g_x\cdot o\in X.
\end{equation}
By definition the operators $\mathcal R$,
$\mathcal R^t$ and $\mathcal R^t \mathcal R$ are
all  $G$-invariant with
respect to the corresponding actions of $G$,
in particular, we have
$$
(\mathcal R^t \mathcal R)f(gx)=
(\mathcal R^t \mathcal R)(f\circ g)(x)$$
for any $g\in G$.

We will write $\mathcal R^t \mathcal R$ as
an integration on $X$. For notational convenience 
we introduce
\begin{equation}
  \label{eq:sinht}
\SH(t)=\prod_{j=1}^r \sh t_j, \,\,\,
\CH(t)=\prod_{j=1}^r \ch t_j, \,\,\,
\,\,\, t=t_1\xi_1+\cdots+ t_r \xi_r \in \fa
\end{equation}
and with some abuse of notation we denote also $|\SH(x)|$ and
$\CH(x)$ the corresponding $K$-invariant
function on $X$.

\begin{lemm+} The operator $\mathcal R^t\mathcal R$
is  given by
$$
\mathcal R^t\mathcal Rf(x)=2^{ra(r^\prime -n)}\int_X
|\SH(\xi )|^{a(r^\prime -n)} f(g_x\cdot \xi)d\mu(\xi), 
\quad x=g_x\cdot o\in X,
\quad f\in C_0^\infty(X).
$$
\end{lemm+}
\begin{proof} We apply the integral
formula (\ref{eq:nt-yo})
to   (\ref{eq:square-radon}),
\begin{equation*}
(\mathcal R^t \mathcal R) f (x)=\int_{K}\int_{K_0}\int_{\fa}
f(g_x kk_0 \exp(t)\cdot o)d\mu_0(t)dkdk_0=
\int_{K}\int_{\fa}
f(g_x  k\exp(t)\cdot o)d\mu_0(t)dk
\end{equation*}
The measure $d\mu_0$ on $\fa$ differs $d\mu$
by a factor,
$$
d\mu_0(t)=2^{ra(r^\prime -n)}|\SH(t)|^{a(r^\prime -n)}d\mu(t)
$$
so that the integration can be written as
on $X$  by the Cartan decomposition 
(\ref{cartan-x}) of $X$,
$$
(\mathcal R^t \mathcal R) f (x)=
2^{ra(r^\prime -n)}\int_{X}
|\SH(\xi)|^{a(r^\prime -n)}
f(g_x  \xi)d\mu(\xi).
$$
\end{proof}

\section{Bernstein-Sato type formula for the function
$|\SH(t)|^\delta$}

\subsection{Invariant differential operators on $X$
}

We recall briefly the relation between invariant
differential operators on $X$ and
Cherednik operators on $\fa$, the latter
being used in the next subsection; the results
are well known, see e.g. \cite{Hkm-Sch-book},
\cite{Opdam-acta},
\cite{Heckman-Opdam-1} and \cite{Heckman-Opdam-2}.

Let $\mathcal Q$ be an invariant differential
operator on $X$. If $f$ is a $K$-invariant
function on $X$ then so is $\mathcal Qf$,
and there is a Weyl group $W$-invariant differential
operator $\text{rad}(\mathcal Q)$ on $\fa$ such that
$(\mathcal Qf)|_{\exp\fa\cdot o}=
\text{rad}(\mathcal Q)f|_{\exp\fa\cdot o}$. 
Let $\{D_j\}$ be the Cherednik operator on $\fa$
for the root system $\Sig(\fg, \fa)$, see below. Then
there is a Weyl group $W$-invariant polynomial $p$ so that
$$
\text{rad}(\mathcal Q)=p(D_1, \cdots, D_r).
$$
In this sense the invariant differential
operator $\mathcal Q$ is uniquely determined
by the polynomial $p$ and vice versa. In particular
the Harish-Chandra homomorphism of $\mathcal Q$
(namely its eigenvalue on  spherical functions)
is just $p$.

\subsection{Bernstein-Sato type formula on $\fa$}

In this subsection we will prove certain Bernstein-Sato
type formula for the function $\SH(t)$ on $\br^r$.
We consider a  root system $\Sig$
of type B, D or  BC as in Section 2
 with general root multiplicity
$a$, $2b$ and $\iota$. The multiplicities
of the root system $\Sig=\Sig(\fg, \fa)$ of  $X$ correspond
to 
\begin{equation}
  \label{eq:b-iota}
(a, 2b, \iota)=\begin{cases} (1, n-2k, 0),&\quad \bbK=\br \\
 (2, 2(n-2k), 1),&\quad \bbK=\bc\\
 (4, 4(n-2k), 3),&\quad \bbK=\mathbb H.\end{cases}
\end{equation}
We assume now $a, b, \iota \in \mathbb C$
and the results in this subsection hold
for this general setting. Let $\rho$ be the  half sum of
positive roots, then
\begin{equation}
  \label{eq:rho}
\rho=\sum_{j=1}^r \rho_j\xi_j^\ast, \quad \rho_j=\iota+ b +\frac a2(r-j).\end{equation}

Let 
\begin{equation*}
\begin{split}
D_j&=\partial_j -a\sum_{i<j}\frac{1}{1-e^{-2(t_i-t_j)}}(1-s_{ij})
+ a\sum_{j<k}\frac{1}{1-e^{-2(t_j-t_k)}}(1-s_{jk})+
\\
&+a\sum_{k\ne j}\frac{1}{1-e^{-2(t_j+t_k)}}(1-\sig_{jk})
+2\iota\frac{1}{1-e^{-4t_j}}(1-\sig_{j})
+2b\frac{1}{1-e^{-2t_j}}(1-\sig_{j})
-\rho_j
\end{split}
\end{equation*}
be the Cherednik operators acting on functions $f(t)$
on $\fa=\br^r$.
Here $s_{ij}, \sig_{ij}, \sig_i$ are the elements in the 
Weyl group,  $s_{ij}=(ij)$ being the permutation
of $\xi_i$ and $\xi_j$, $\sig_{ij}$ the signed permutation,
$\sig_{ij}(\xi_i)=-\xi_j, \sig_{ij}(\xi_j)=-\xi_i$,
and $\sig_{i}$ the reflection $\sig_i(\xi_i)=-\xi_i$,
all of these mapping $\xi_k\to \xi_k$ for $k\ne i, j$; 
see \cite{Opdam-acta} in general case
and \cite{gz-imrn} for root systems of Type BC.
(Note that we are using some different
convention of root system here from that  in \cite{Opdam-acta}.
The roots here are half of those  in \cite{Opdam-acta}
whereas the multiplicities are twice of those there.)

\begin{theo+}\label{BSF} Let 
$m_\delta$ and $\mathcal M_{\delta}$ 
be the following constant and respectively $W$-invariant polynomial
of the Cherednik operators $\{D_j\}$,
$$
m_\delta=\prod_{j=1}^r \big(\delta +a(j-1) \big)
\big( \delta-1 +\iota +2b +a(r-j)\big), 
\quad
\mathcal M_\delta=\prod_{j=1}^r
\big(D_j^2-(\delta+\rho_1)^2\big),
\quad \delta\in \mathbb C.
$$
Then the following Bernstein-Sato type identity holds
\begin{equation}
  \label{BS-sinh}
\mathcal M_\delta | \SH(t)|^{\delta}
=m_\delta
|\SH(t)|^{\delta-2}, \quad t=(t_1, \cdots, t_r), \, t_j\ne 0 \quad\quad \forall j
\end{equation}
\end{theo+}

\begin{proof} Following the proof of a
Bernstein-Sato type formula 
(see (\ref{BS-cosh}) below) for 
$\CH^{\delta}$ 
in \cite{gz-imrn} we factorize the operator $\mathcal M_\delta$
as 
$$
\mathcal M_\delta=\prod_{j=1}^r \big(D_j -(\delta+\rho_1)\big)
\prod_{j=1}^r \big(D_j+\delta+\rho_1\big),
$$
and consider  successively the actions
of $(D_j+\delta+\rho_1\big)$ from $j=1$ to $j=r$
and the actions of
$(D_j-(\delta+\rho_1)\big)$ backward from $j=r$ to $j=1$.
\begin{lemm+} For $j=1, \cdots, r$ 
 we have
\begin{equation*}
  \begin{split}
&\quad 
\big(\prod_{k=1}^j (D_k+\delta+\rho_1)\big)|\SH(t)|^{\delta} \\
&=\big(\prod_{k=1}^j 
(\delta +a(k-1))\big)
|\SH(t)|^{\delta} \prod_{k=1}^{j} 
(1+\coth t_k),
  \end{split}
\end{equation*}
and
\begin{equation*}
  \begin{split}
&\quad 
\big(\prod_{k=j}^r (D_k-(\delta+\rho(\xi_1))\big)
\big(|\SH(t)|^{\delta}\prod_{i=1}^r (1+\coth t_i)\big )
\\
&=\big(\prod_{k=j}^r(\delta-1 +(r-k)a +\iota +2b)\big)
|\SH(t)|^{\delta} \prod_{i=1}^r (1+\coth t_i)
 \prod_{k=j}^r(\coth t_k-1),  \end{split}
\end{equation*}
where  $t=(t_1, \cdots, t_r)$, $t_k>0 \forall k$ 
\end{lemm+}

\begin{proof} 
 The first  formula
can be proved by almost
the identical  computations as in the proof of Lemmas 2.2
\cite{gz-imrn} 
changing $\ch$ to $\sh$, $\tanh $ to $\coth$,
producing also the same factor
$\prod_{k=1}^j 
(\delta +a(k-1))$.
So is the second formula as the proof of Lemma 2.3, loc. cit.,
however a different factor
will appear; we give
a sketch of the proof for $j=r$
showing the difference. 
 The rest
is done by the backward induction
and some similar  yet long computations.
For brevity
we denote temporarily
$P=\prod_{j=1}^r(1+\coth t_j)$
and $\hat P_{I} =\prod_{l\notin I}(1+\coth t_l)$
for any index subset $I$.
We compute $D_r\big( |\SH(t)|^{\delta}P\big)$
and write the result in terms of $ |\SH(t)|^{\delta}P$
and $|\SH(t)|^{\delta}P(\coth t_r -1)$.
Using 
$$
\frac{d}{ds}|\sh(s)|^{\delta}=\delta |\sh(s)|^{\delta} (\coth s -1)+
\delta |\sh(s)|^{\delta}, \quad
\frac{d}{ds}\coth (s)=-(\coth s +1)(\coth s -1)
$$ 
and that $P$ is invariant under permutations,
we have
\begin{equation*}
\begin{split}
&\quad\,D_r\big( |\SH(t)|^{\delta}P\big)\\
&=\delta |\SH(t)|^{\delta} P (\coth t_r -1)
+ \delta |\SH(t)|^{\delta} P 
-|\SH(t)|^{\delta} \hat P_r (\coth t_r +1) (\coth t_r -1)
\\
& \quad + 2a |\SH(t)|^{\delta} 
\sum_{i<r}\hat P_{i,r}\frac{e^{t_i+t_r}}{e^{t_i +t_r}-
e^{-(t_i +t_r)}} (\coth t_r+\coth t_i)
\\
&\quad +2\iota |\SH(t)|^{\delta} \hat P_r
\frac{e^{2t_r}}{e^{2t_r}-e^{-2t_r}} 
2\coth t_r 
+2b |\SH(t)|^{\delta} 
\hat P_r
\frac{e^{t_r}}{e^{t_r}-e^{-t_r}} 
2\coth t_r 
-\rho_r |\SH(t)|^{\delta}P.
\end{split}
\end{equation*}
The third term is
$$-|\SH(t)|^{\delta} \hat P_r (\coth t_r +1)(\coth t_r -1)
=-|\SH(t)|^{\delta} P(\coth t_r -1).$$
Each term in the summation $\sum_{i<r}$ is
\begin{equation*}
\begin{split}
&\quad\,\hat P_{i,r}
\frac
{(\ch t_i +\sh t_i)(\ch t_r +\sh t_r)}
{2\sh(t_i +t_r)} 
(\coth t_i +\coth t_r) \\
&=\frac 12 \hat P_{i,r}
\frac
{(\ch t_i +\sh t_i)(\ch t_r +\sh t_r)}
{\sh(t_i +t_r)} \frac{\sh(t_i +t_r)}
{\sh t_i \sh t_r}\\
&=\frac 12\hat P_{i,r}(\coth t_r +1)(\coth t_i +1)=\frac 12 P.
\end{split}
\end{equation*}
In the next summand we notice that
\begin{equation*}
\begin{split}
\hat P_r\frac{e^{2t_r}}{e^{2t_r}-e^{-2t_r}} 
\coth t_r 
&=\hat P_r\frac{\sh(2t_r) +\ch(2t_r)}{2\sh(2t_r)}
\coth t_r \\
&=\hat P_r\frac{\sh(2t_r) +\ch(2t_r)}{2\sh(2t_r)}
\coth t_r 
=\frac 14 \hat P_r(\coth^2 t_r -1 +2(\coth t_r +1))\\
&=\frac 14 (P (\coth t_r -1) +2P).
\end{split}
\end{equation*}
Similarly, in the  next  term we have
$$
P_r
\frac{e^{t_r}}{e^{t_r}-e^{-t_r}} \coth t_r 
=\frac 12 (P (\coth t_r -1)+ P)
$$
and then
\begin{equation*}
\begin{split}
&\quad\,\big(D_r -(\delta +\rho_1)\big)
\big( |\SH(t)|^{\delta}P\big)\\
&=\big(\delta + 2\iota +2b +(r-1)a -\rho_r -(\delta +\rho_1)\big)
|\SH(t)|^{\delta}P \\
&\qquad
+(\delta -1 +\iota +2b)|\SH(t)|^{\delta}P (\coth t_r -1)
\end{split}
\end{equation*}
but the first term is vanishing and this
proves the identity for $j=r$.
\end{proof}

Our result  follows then from the  first formula for $j=r$ and the second
for $j=1$,
\begin{equation*}
  \begin{split}
\mathcal M_\delta |\SH(t)|^{\delta}
&=\prod_{k=1}^r \big(D_k-(\delta+\rho_1)\big)
\prod_{k =1}^r \big(D_k+\delta+\rho_1\big)
|\SH(t)|^{\delta}\\
&=\big(\prod_{k=1}^r
(\delta +a(k-1)) \big)
\prod_{k=1}^r \big(D_k-(\delta+\rho_1)\big)
\big(|\SH(t)|^{\delta} \prod_{i=1}^r (\coth t_i +1)\big )
\\
&=\big(\prod_{k=1}^r(\delta +a(k-1)) \big) 
\big(\prod_{k= 1}^r(\delta-1 +(r-k)a +\iota +2b)\big)
\\
&\quad \times |\SH(t)|^{\delta}\prod_{i=1}^r (\coth t_i+1)
\prod_{k=1}^r (\coth t_k-1)
\\
&=m_\delta |\SH(t)|^{\delta-2}.
\end{split}
\end{equation*}
\end{proof}

It is worthwhile to make some remarks on the theorem.
\begin{rema+} \label{sh-ch-rem}
In \cite{gz-imrn} the following identity
is proved
\begin{equation}
  \label{BS-cosh}
\mathcal M_\delta \CH(t)^{\delta}
=(-1)^r\prod_{j=1}^r \big(\delta +a(j-1) \big)
\big( \delta-1 +\iota+a(r-j) \big)
\CH(t)^{\delta-2},
\end{equation}
which is of fundamental importance
in Berezin transform \cite{gkz-manu-mat}
and in spherical transforms
of Jacobi type functions \cite{gz-imrn}.
The pair of equalities 
  (\ref{BS-sinh}) and 
(\ref{BS-cosh}) is to be compared with the following
pair of simple identities (for the flat rank one case with all root
multiplicities being zero)
\begin{equation}
  \label{eq:BS-flat}
(\frac{d^2}{dt^2}-\delta^2) |\sh(t)|^{\delta}
=\delta(\delta -1)|\sh(t)|^{\delta-2}, \quad (\frac{d^2}{dt^2}-\delta^2) \ch(t)^{\delta}
=-\delta(\delta -1)\ch(t)^{\delta-2}.
\end{equation}
Note that in the flat case for $\delta$ being an even integer
it is possible to derive one equality in the pair
 from the another by using analytic continuation
in $t$,  the fact that $\sh(t+\frac{\pi}2 i)=i \ch t$
and that the operator $\frac{d^2}{dt^2}$
is translation invariant. However 
in the non-flat case the two equalities 
  (\ref{BS-cosh})
 and  (\ref{BS-sinh})
  are considerably different,
even though the methods of the proofs are
similar.
Firstly the Cherednik operators are not
translation invariant
and it seems to the author that
 intertwining relations between
the Cherednik operators and translations
are not fully understood (see \cite{Cherednik-imrn97}
for some related questions); secondly
the constants in the right hand
of   (\ref{eq:BS-flat}) differ
just by a sign change and it is not
the case for   (\ref{BS-cosh})
 and  (\ref{BS-sinh}) when  $b\ne 0$.
Also, the equality 
  (\ref{BS-sinh}) holds true for $t$ with $t_j\ne 0$
as the function $|\SH(t)|^{\delta}$ has singularity
on the hypersurfaces $t_j=0, j=1, \cdots, r$ for general $\delta\in \bc$
and is not differentiable, whereas $\CH(t)^{\delta}$
is a real  analytic function.

When the root 
 multiplicities correspond to
the Grassmannian manifold
in $\bbK^n$
and for certain specific integers $\delta$
(depending on the domain)
it should be possible to  derive
some Bernstein-Sato type formula
for the function
$|(\prod \sin t_j)|^{\delta}$ on the compact torus
$i\br^r/i\pi \mathbb Z^r$
by combining the result
of Grinberg \cite{Grinberg-grs}
on Radon transform on $X^c$
and the result of Opdam \cite{Opdam-acta}
on eigenvalues of the Cherednik operators,
and it is interesting to carry out the details;
see also \cite{Kakehi-jaf99} and \cite{Gonzales-Kakehi-1}.
\end{rema+}

\subsection{The meromorphic
continuation of the
distribution $|\SH(t)|^{\delta}$ 
on root systems with positive
multiplicities}

 We assume $\Sig$ is as
in the previous subsection  a root
system with general non-negative multiplicities
$a\ge 0$, $b\ge 0$  and $\iota\ge 0$.

We notice first that, if $f\in  C_0^\infty(\br^r)^W$ namely Weyl group invariant,
 the integral
$$
\int_{\fa}|\SH(t)|^{\delta} f(t)d\mu(t)
$$
is absolute convergent if $\Re(\delta)>-1-l-2b$, and
it defines a distribution analytic in $\delta$.
We consider now the integral
\begin{equation}
  \label{eq:zeta-del}
\zeta_{\delta}(f)
=\frac{1}{z_\delta}
\int_{\fa}|\SH(t)|^{\delta} f(t)d\mu(t)
\end{equation}
where 
\begin{equation}
{z_\delta}={\Gamma_{a}(\frac {\delta}2 +\frac a2(r-1)+1 )
\Gamma_{a}
(\frac {\delta+\iota +2b }2 +\frac a2(r-1) +1)},
\end{equation}
and
\begin{equation}
\label{gg++}
\Gamma_{a}(\lam)
=\prod_{j=1}^r\Gamma(\lam-\frac a2(j-1))
\end{equation}
is the Gindikin Gamma function.
\begin{theo+} \label{an-cont}
 The distribution $\zeta_\delta$
on $C_0^\infty(\br^r)^W$ for $\Re\delta
>-1-\iota -2b$ has analytic continuation
to  $\delta \in \mathbb C$ and 
$$
\mathcal M_{\delta} \zeta_{\delta} =2^{2r}\zeta_{\delta -2},
$$
where $\mathcal M_{\delta}$ is the differential operator in Theorem 3.1.
\end{theo+} 

\begin{proof} We rewrite  the factor $m_\delta$ in Theorem \ref{BSF} using
$\zeta_{\delta -2}$.
Changing $j-1$ to $r-j$ we see that
$$
m_\delta=2^{2r} \prod_{j=1}^r \big(\frac{\delta}2
+\frac a2 (r-1)-\frac a2 (j-1) \big)
\prod_{j=1}^r \big( \frac{\delta-1+\iota +2b}2 +\frac a2 (r-1)-\frac a2 (j-1) \big),
$$
and each product is a quotient of the Gindikin Gamma function
\begin{equation}
\label{guo-gamm-1}
 \prod_{j=1}^r \big(\frac{\delta}2
+\frac a2 (r-1)-\frac a2 (j-1) \big)=\frac{\Gamma_a( \frac{\delta}2
+\frac a2 (r-1) +1)}
{\Gamma_a( \frac{\delta}2
+\frac a2 (r-1))},
\end{equation}
\begin{equation}
\label{guo-gamm-2}
\prod_{j=1}^r \big( \frac{\delta-1+\iota +2b}2 +\frac a2 (r-1)-\frac a2 (j-1) \big)
=\frac{\Gamma_a(  \frac{\delta-1+\iota +2b}2 +\frac a2 (r-1) +1)}
{\Gamma_a(\frac{\delta-1+\iota +2b}2 +\frac a2 (r-1))}.
\end{equation}
Namely $m_\delta=2^{2r}
=2^{2r}\frac{z_\delta}{z_{\delta-2}}$ and the Bernstein-Sato formula
  (\ref{BS-sinh}) is then
$$
\frac{1}{{z_\delta}}
\mathcal M_{\delta} 
|\SH|^{\delta}
=2^{2r}
\frac{1}{z_{\delta-2}}
|\SH|^{\delta-2}.
$$
The function $|\SH|^{\delta}$ 
and $\mathcal M_\delta |\SH|^{\delta}$ 
for sufficient large
$\delta$ are  well-defined distributions and
are analytic in $\delta$, and
so is the constant
$\frac{1}{{z_\delta}}$.
 By partial integration
we have
\begin{equation*}
\begin{split}
\zeta_{\delta}(f)&=\frac{1}{{z_\delta}}\int_{\mathbb R^r}|\SH(t)|^{\delta}\mathcal M_\delta f (t) d\mu(t)\\
&=\frac{1}{{z_\delta}}\int_{\mathbb R^r} \mathcal M_\delta |\SH(t)|^{\delta}\mathcal M f (t) d\mu(t)\\
&=2^{2r}\zeta_{\delta-2}(f)
\end{split}
\end{equation*}
see the proof of \cite[Lemma 7.8]{Opdam-acta}.
This proves the analytic continuation.
\end{proof}

See also  Remark \ref{last-rem} for a possible refinement of the above theorem.

\section{The Dirac distribution
as meromorphic continuation of the distribution 
$|\SH(t)|^{\delta}$ on symmetric domains.}

\subsection{G\aa{}rding-Gindikin distribution
$(x_1 \cdots x_r)_+^{\lam-\frac a2(r-1)-1}$ on $\mathbb R^r$}

In this subsection we reformulate
the result of Gindikin \cite{FK-book}
on Riesz type integrals (also called
G\aa{}rding-Gindikin distributions)
on the Jordan algebra
of  real  symmetric,
complex hermitian and quanternionic hermitian
$r\times r$ matrices) in the polar coordinates with the test function
being $U(r, \bbK)$-invariant smooth functions $f$
of compact support.
(In \cite{FK-book} the more general case of
Schwartz functions is considered.)  We will identify such
functions $f$ as a function on $\mathbb R^r$.
 We refer to \cite{FK-book}
for an account of symmetric cones; see also
 \cite{Ourn-Rubin-amsm}.
Let
$$(x_1 \cdots x_r)_+
=(x_1)_+ \cdots (x_r)_+
$$
where  $x_+$ on $\br$
is the function $x_+=x$ if $x>0$ and $x_+=0$
is $x\le 0$.
Consider the G\aa{}rding-Gindikin integral
\begin{equation}
  \label{eq:gg-int}
  \begin{split}
\mathcal G_{\lam}(f)&=\frac 1{r!} \frac 1{\Gamma_a(\lam)}\int_{\mathbb R^r}
(x_1 \cdots x_r)_+^{\lam-\frac a2(r-1)-1} f(x_1, \cdots , x_r)
\prod_{i<j}|x_i-x_j|^{a}
dx\\
&= \frac 1{r!}\frac 1{\Gamma_a(\lam)}\int_{\mathbb R_+^r}
(x_1 \cdots x_r)^{\lam-\frac a2(r-1)-1} f(x_1, \cdots , x_r)
\prod_{i<j}|x_i-x_j|^{a}
dx
  \end{split}
\end{equation}
where $f\in C_0(\br^r)$,
is symmetric in $x_1, \cdots, x_r$,
and
$\Gamma_a(\lam)
$
is the Gindikin Gamma function
 (\ref{gg++}), and $a=1, 2, 4$.

The following lemma follows
immediately from \cite[Theorem VII.2.2]{FK-book}
and from (a special case of) Chevalley's  theorem
that the restriction to the diagonal of $U(r, \bbK)$-invariant
smooth functions of compact support on the space of all
symmetric matrices
to the diagonal matrices is an isomorphism (see
\cite[Chapter II, Section 5]{He2} for some general statements).

\begin{lemm+} \label{an-GG}
Suppose $f\in C_0^\infty(\br^r)$  is symmetric.
The distribution $\mathcal G_{\lam}$ has analytic continuation
to $\lam\in \mathbb C$ and
$$
\mathcal G_{0}=c_0 \mathcal I_0.
$$
where
\begin{equation}
c_0=\prod_{1\le i<i\le r}\frac{\Gamma(\frac a2(j-i+1))}
{\Gamma
(\frac a2(j-i))}
\end{equation}
and $\mathcal I_0$ is the Dirac distribution at $x=0$.
\end{lemm+}

\subsection{Analytic
continuation of $|\SH(t)|^{\delta}$.
The case of symmetric domains}

We let $\Sig$ be as in Section 2
the root system of the symmetric
domain $X$ with the multiplicity
$(a, 2b, \iota)$ given as in  (\ref{eq:b-iota}).
\begin{prop+}\label{prop-dirac} Suppose that
the root multiplicities satisfy
 $1+\iota +2b>r(a-1)$. 
Let $\delta_0\in \br$  be such that
$$
-1-\iota-2b <\delta_0 < -r(a-1)
$$
and that
$$
l=\frac{\delta_0 +\iota +2b -1 +a(r-1)}2 +1
$$
is a positive integer.
Then for any $f\in C_0^\infty(\fa)^W$,
the following formula holds
$$
\left(\prod_{k=0}^{l-1}\mathcal M_{\delta_0-2k}\right)
|\SH(t)|^{\delta_0} =c_1\mathcal I_0
$$
in the sense of distribution on 
$C_0^\infty(\fa)^W$ (defined by the integration
against $d\mu(t)$),
where $\mathcal I_0$ is the Dirac distribution at
$t=0$,
$$
c_1=2^{r(2b+2\iota)+lr } r!
\Gamma(\frac 12(\delta_0 +\iota +2b -1)+\frac a2(r-1)+1) 
(\prod_{k=0}^{l-1}\prod_{j=1}^r(\delta_0 -2k +a(j-1))
c_0 
$$
and $c_1\ne 0$.
\end{prop+}
\begin{proof} Let  $\Re \beta>2l-\iota-2b+1$ be sufficiently large.
We compute
\begin{equation*}
\begin{split}
S_\beta(f):&=\int_{\mathbb R^r}|\SH(t)|^{\beta}
\left(\prod_{k=0}^{l-1}\mathcal 
M_{\beta -2k}\right) f (t) d\mu(t)\\
&=\int_{\mathbb R^r}\big( \left(\prod_{k=0}^{l-1}
\mathcal M_{\beta -2k}\right) |\SH(t)|^{\beta}\big) f (t) d\mu(t)\\
&= (\prod_{k=0}^{l-1}m_{\beta-2k}) \int_{\br^r}|\SH(t)|^{\beta-2l}
 f (t) d\mu(t)\\
&= (\prod_{k=0}^{l-1}
m_{\beta-2k})2^r \int_{\br_+^r}|\SH(t)|^{\beta-2l}
 f (t) d\mu(t).
\end{split}
\end{equation*}
Change variables $x_j=\sh t_j^2$. Then
\begin{equation}\label{S-GG}
\begin{split}
S_\beta(f)&= 2^{r(2b+2\iota) }(\prod_{k=0}^{l-1}m_{\beta-2k})
 \int_{\br_+^r}
(x_1\cdots x_r)^{\frac 12(\beta-2l +\iota +2b -1)}
\\
&\times \prod_{j=1}^r(1+x_j)^{\frac{1}2(\iota -1)}
f(x_1, \cdots, x_r) \prod_{i<j}(x_i -x_j)^a
dx_1 \cdots dx_r\\
&= 2^{r(2b+2\iota) } r! 
(\prod_{k=0}^{l-1}m_{\delta-2k})
\Gamma_a(
\frac 12(\beta-2l +\iota +2b -1)+\frac a2(r-1)+1)\\
&\times \mathcal G_{\frac 12(\beta -2l+\iota +2b -1)+\frac a2(r-1)+1}(F)
\end{split}
\end{equation}
in terms of the G\aa{}rding-Gindikin distribution,
where 
$$
F(x_1, \cdots, x_r)= \prod_{j=1}^r(1+x_j)^{\frac{1}2(\iota -1)}
f(x_1, \cdots, x_r)
$$ 
for $x\in \mathbb R_+^r$
and is extended to any symmetric function 
in $C^\infty_0(\mathbb R^r)$.
Using  (\ref{guo-gamm-2}) we see that the
in the previous formula
the constant in front of the G\aa{}rding-Gindikin distribution is
(disregarding  $2^{r(2b+2\iota) } r!$ )
\begin{equation*}
\begin{split}
&\quad\,(\prod_{k=0}^{l-1}m_{\delta-2k})
\Gamma_a(\frac 12(\beta-2l +\iota +2b -1)+\frac a2(r-1)+1)\\
&=2^{lr} \Gamma_a(\frac 12(\beta +\iota +2b -1)+\frac a2(r-1)+1)
(\prod_{k=0}^{l-1}\prod_{j=1}^r(\delta -2k +a(j-1))
\end{split}
\end{equation*}
is analytic for all $\beta\in \bc$
such that
\begin{equation}
  \label{eq:con-beta}
\frac 12(\Re\beta +\iota +2b -1)+1>0  
\end{equation}
Thus the above
formula (\ref{S-GG}) hold for all such $\beta$.
Our results follows by taking
 $\beta=\delta_0$ and applying Lemma \ref{an-GG}. The assumption on
the root multiplicities, on $\delta_0$ and $l$ guarantees
that the constant $c_1\ne 0$.
\end{proof}

\section{Right inverse of $\mathcal R^t\mathcal R$.
Left inverse of $\mathcal R$ in the case of rank one}

\subsection{Right Inverse of the transform $\mathcal R^t R$
}

\begin{theo+}\label{main-thm}
Let $X$ be the symmetric domain $X=U(n-r, r; \bbK)/U(n-r;\bbK)
\times U(r,\bbK)$, and
 $Y=U(n-r, r; \bbK)/U(n-r^\prime; \bbK)\times U(r^\prime; \bbK)$
be the set of $G$-translates of the totally geodesic submanifold
$y_0$.
Suppose that the rank condition (\ref{ran-cd}) is
satisfied, and that $r^\prime -r$ is even
when $\bbK=\br$.
Let $\mathcal M$ be the invariant differential
operator on $X$ so that
$$
\rad \mathcal M
=\prod_{k=0}^{l-1}\mathcal M_{a(r\prime -n)-2k},
\quad l=\frac a2(r^\prime-r)
$$
Then
\begin{equation}
  \label{eq:rgt-ivs}
\mathcal R^t \mathcal R \mathcal M  f(x)=cf(x), 
\quad f\in C_o^\infty(X), \qquad x\in X.
\end{equation}
where
$c=c_1$ as given in Propostition 4.2 with $\iota=a-1, 2b=a(n-2r)$
and $\delta_0=a(r^\prime-n)$.
 \end{theo+}

 \begin{proof}Suppose first that $f\in C_0^\infty(X)$ is $K$-invariant.
Then by Lemma 2.2,
$$
(\mathcal R^t\mathcal R \mathcal M f)(o)=
2^{ra(r^\prime -n)}\int_{\fa}|\SH(t)|^{\delta_0} \rad(\mathcal M) f(\exp(t)\cdot o)
d\mu(t)
$$
with $\delta_0=a(r^\prime -n)$. 
Applying Proposition 4.2 we find
$$
(\mathcal R^t\mathcal R \mathcal M f)(o)=
cf(o).
$$

Now for any $f\in C_0^\infty(X)$ and for a fixed  $x=g_x\cdot o$ we let
$$
F(\xi)=\int_K f(g_x k\xi)dk,
$$
which is a $K$-invariant function in $C_0^\infty(X)$.
The inversion formula for $F$ at $\xi=o$ reads as follows
$$
 (R^t R \mathcal M  F)(o)
=c F(o)=cf(x)
$$
However by invariance of of $R^t R$
and  $\mathcal M$ under $G$ we
have
\begin{equation*}
  \begin{split}
   (R^t R \mathcal M  F)(o)&= \int_K 
\big((R^t R \mathcal M)f(g_x k\xi)\big)\big|_{\xi=o}dk
=\int_K (R^t R \mathcal Mf)(g_x k\xi)\big|_{\xi=o}dk\\
&=\int_K (R^t R \mathcal Mf)(g_x ko)dk
=(R^t R \mathcal Mf)(x),
  \end{split}
\end{equation*}
completing the proof.
 \end{proof}

We note that the above theorem 
is closely related to the surjectivity
of the invariant differential operator $\mathcal M$.
We let $\mathcal M$ act
on the formula (\ref{eq:rgt-ivs}), 
$$
\mathcal M\mathcal R^t \mathcal R u(x)=c u(x),
$$
where $u=\mathcal Mf$.
Thus we have an inversion formula for functions
$u\in C_0^\infty(X)$ which are in the image 
of $\mathcal M$ on $C_0^\infty(X)$. Generally
speaking however, the image of
an invariant differential operator
on $C_0^\infty(X)$ is not onto. Recall,
 by a theorem Helgason\cite{He-annmath73}, that
for any $u\in C_0^\infty(X)$ (and
for any invariant differential operator, in particular
for $\mathcal M$)
there is an $f\in C^\infty(X)$ such that $\mathcal Mf=u$; however
the function $f$ is not necessarily in $C_0^\infty(X)$.

\subsection{The case of rank one}

When  $r=1$ the matrix space
$M_{1, n-1}(\mathbb K)=\bbK^{n-1}$,
and the  domain $X$ is the hyperbolic
ball in $\bbK^{n-1}$.
The subdomain
$y_0$ in question is  the unit ball in $0\oplus \bbK^{r^\prime-1}$.
We consider the Radon transform
over the set $Y=G\cdot y_0$ of $(r^\prime -1)$ dimensional
totally geodesic submanifolds, for any $1<r^\prime \le n-1$.
We
can also include the exceptional
domain
into our consideration, with $X$ being the unit ball
 in $\bbK^2$ where
$\bbK=\mathbb O$, real Cayley numbers,  and 
$y_0$  the unit ball in $0\oplus \mathbb O=\mathbb R^8$, $y_0=SO_0(8, 1)/SO(8)$.
As  symmetric spaces $X=G/K$, $y_0=G_0/K_0$
with 
$$
(\fg, \fk)=(\ff_{4(-20}, \mathfrak{so}(9)),\quad
(\fg_0, \fk_0)=(\mathfrak{so}(8, 1), \mathfrak{so}(8)).
$$
The parameters $a, r, r^\prime$
are given by
$$
(a, r, r^\prime)=(8, 1, 2).
$$

We specialize now  the Bernstein-Sato type formula
in the rank one cases, with the 
 root multiplicities there given by
$$
(2b, \iota)=(a(n-2), a-1)
$$
with the convention that $\mathbb K^{n-1}=\mathbb O^2$
and $a(n-2)=8$ when $a=8$ i.e.
$\mathbb K=\mathbb O$. 
The half sum of positive
roots is $\rho=b +\iota=\frac{a}{2}n -1$.

We normalize the Riemannian metric
on $X=G/K$ so that the vector $\xi=\xi_1\in \fp$,
viewed as a tangent vector of $X$ at $o$, has norm $1$;
we use the same Riemannian measure $d\mu(x)$ as
above.
The radial part $\rad(\mathcal L)$ of the Laplace-Beltrami operator
$\mathcal L$ is  then
\begin{equation}
  \label{eq:norm-lap}
\rad(\mathcal L)= D^2-\rho^2  
\end{equation}
where $D$ is Cherednik operator; this follows
directly by computing $D^2$ and by using 
the formula in \cite[Chapter II]{He2}
for the  radial part
of Laplace-Beltrami operator\footnote{The Riemann metric used \cite{He2}
is defined by the Killing form, so that
the tangent vector $\xi$ has squared norm 
$\tr_{\fg}(\ad\xi\ad \xi)=4(an-2)$, so that our Laplace-Beltrami
operator here is that in \cite{He2} multiplied by
$4(an-2)$.}
 on the symmetric
space $X$.

We let
\begin{equation}
\label{helmhotz}
\lam_j=\big(a(n-1)-a(r^\prime -1) +2(j-1)\big )\big(a(r^\prime -1) +a -2j\big), \quad j=1, 2. \cdots.
\end{equation}

\begin{theo+}\label{thm-rk1}
 Let $\mathbb K=\mathbb R, \bc, \mathbb H, \mathbb O$
be the ring of   real, complex, quaternionic and Cayley
numbers.
Let $X=G/K$ be the unit ball in
$\mathbb K^{n-1}$ 
and $Y=G\cdot y_0$ the space of $r^\prime -1$-dimensional
geodesic submanifolds of $X$. Suppose $r^\prime-1$
 is even in the real case $\bbK=\br$ 
(namely $y_0$ is an even $(r^\prime-1)$-dimensional submanifold).
Let  
$$
l=\frac a2(r^\prime-1)
$$ 
Then the Radon  transform $\phi=\mathcal Rf$,
 for $f\in C_0^\infty(X)$
 is inverted by
$$
\mathcal M (\mathcal R^t \phi) =c f
$$
where
$$
\mathcal M=\prod_{j=1}^{l}
(\mathcal L + \lam_j)
$$
$\mathcal L$ is the Laplace-Beltrami operator
(normalized as above).
 \end{theo+}

\begin{proof} We
use Theorem \ref{BSF}. The  operator $\mathcal M_\delta$ there
is, by  (\ref{eq:norm-lap}),
 $$
\mathcal M_\delta= D^2-(\delta+\rho)^2
=D^2-\rho^2 +\rho^2-(\delta+\rho)^2
=\rad(\mathcal L) +(-\delta)(2\rho+\delta)
$$
We have,
\begin{equation*}
\big(\mathcal L+(-\delta)(2\rho+\delta)\big)
|\SH(x)|^{\delta}\\
=\delta(\delta-1+\iota+2b )  |\SH(x)|^{\delta-2}, \quad x\ne o
\end{equation*}
which can also be proved
independently by straightforward computations. We recall
that $|\SH(x)|$ is extended to a $K$-invariant function on $X$.
With this formula it is now easy to study the analytic continuation
of the distributions by elementary computations
independent of the previous sections. 
In particular, by Proposition \ref{prop-dirac}
 we have 
\begin{equation}\label{M-sh}
\mathcal M |\SH(x)|^{a(r^\prime -n)} 
=c_3 \mathcal I_o
\end{equation}
where $\mathcal I_o$ is the Dirac distribution at the origin
$o$. 

The rest is essentially the same as
\cite[Chapter I, Section 4]{He2} (with some different formulation).
Let the operator $\mathcal M$ act on the formula
$\mathcal R^t\mathcal R f(x)$ in Lemma 2.2
identified as a function of $g\in G$
\begin{equation*}
  \begin{split}
\mathcal M_x \mathcal R^t\mathcal Rf(x)
&=\mathcal M_g \mathcal R^t\mathcal Rf(g)\\
&=2^{(2b^\prime -2b)}\int_X
|\SH( \xi)|^{a(r^\prime -n)} 
\mathcal M_g f(g\xi)d\mu(\xi)\\
&=2^{(2b^\prime -2b)}\int_X
|\SH( \xi)|^{a(r^\prime -n)} 
(\mathcal M_\xi f)(g\xi)d\mu(\xi)\\
&=2^{(2b^\prime -2b)}\int_X
|\SH( \xi)|^{a(r^\prime -n)} 
(\mathcal M_\xi f)(g_x \xi)d\mu(\xi)    
  \end{split}
\end{equation*}
by the   bi-invariance of the Laplace-Beltrami operator.
The rest follows from
(\ref{M-sh}), and the constant $c$ is the 
same as in Theorem \ref{main-thm}.
\end{proof}

It might be illuminating to write
the operator $\mathcal M$
as
$$
\boxed{\mathcal M=\prod_{j=1}^{\frac 12 d_0}
\big(\mathcal L + (d-d_0 +2(j-1))
(d_0 +a -2j)\big)}
$$
where $d=\dim(X)$ and $d_0=\dim(y_0)$.
The case when  $r^\prime =n-1$ 
was proved proved earlier by
 Helgason \cite{He2} by  case by case computations
(our $d=\dim(X)$ and $r^\prime-1$ corresponds to
$n$ and respectively $n-1$ there).
The roots of the operator $\mathcal M$
as a polynomial of the Laplace-Beltrami
operator in the complex case
appear also in the context
of Helmholtz operators \cite{Schim-Schli-acta};
in particular using our Theorem \ref{BSF}
we can get a more precise form 
of the formula (29) there, but we will go into the details.

Finally we make some remarks on
certain  interesting open questions.

\begin{rema+} Let $\Sig$ be a root system with non-negative
multiplicities.
We introduce the Schwartz space $\mathcal S_\Sig^0(\fa)$
on the root system $\Sig $ on 
$\fa$: A Weyl group invariant smooth function $f$ on $\fa$ is in $\mathcal S^0_{\Sig}(\fa)$
if for any $W$-group invariant
polynomial $p$,
\begin{equation}
\label{sch-con}
\sup_{t\in \fa} e^{-N \rho(t))} |p(D_1, \cdots, D_r) f(t)|
<\infty
\end{equation}
for all $N=0, 1, \cdots$. It is easy to see that the function
 $|\SH(t)|^{\delta}$ can be defined as
a distribution
on $\mathcal S_\Sig^0(\fa)$ and
Theorem 3.4 
is still true in this sense.

When the root system corresponds to a Riemannian symmetric space
$X=G/K$ this is then the Schwartz $L^0$-space $\mathcal S^0(G)$
of $K$ bi-invariant
functions on $G$ \cite{GV}. 
Recall \cite{GV} that the Schwartz $L^p$-space
on $X=G/K$ is the
space of $K$-invariant smooth functions
on $X$ such that
for any invariant differential operator $\mathcal L$,
and any $M>0$,
$$
\sup_{t\in \fa} e^{-(2/p-2) \rho(t)}  
(1+|t|^2)^{-M} \phi_0(\exp (t)\cdot o)
|\mathcal L f(\exp (t)\cdot o)| <\infty.
$$
When $p=0$ that condition can be simplified as
$$
\sup_{t\in \fa} e^{-N \rho(t)} |\mathcal L f(\exp (t)\cdot o)|
<\infty.
$$
To see this one has to use
\cite[Theorem 4.6.6]{GV} there, see also \cite[Proposition 2.2.12]{Anker-Ji}.  
The last condition is clearly equivalent to
(\ref{sch-con}). 
\end{rema+}

\begin{rema+}
Note that the question
of finding an inversion formula
for the Radon transform on higher rank matrix domains is
still open. It seems that
the approach in \cite{Grinberg-Rubin}
and \cite{gz-radon} for compact
symmetric domains combined with 
the results in this paper could lead
to some  interesting results.
As is mentioned in Remark \ref{sh-ch-rem} 
the functions $|\CH(x)|^\delta$ for certain $\delta$ 
is the kernel of the Berezin transform which commutes
with invariant differential operators (due
to an explicit formula for $|\CH(x)|^\delta$ involving
certain Bergman kernels)
whereas it is
not known if it is true for the operator
 with the kernel $|\SH(x)|^\delta$ (and we think
it is not).
It would be interesting to understand the precise
relations between the two operators.
\end{rema+}

\begin{rema+}\label{last-rem} The reader may have noticed
  Proposition \ref{prop-dirac} is proved by using
Lemma  \ref{an-GG}, which is only
proved (to the knowledge of the author) for
root systems of symmetric cones. We conjecture
that  Lemma  \ref{an-GG}
 is true for general nonnegative multiplicity $a$
and for general
non-symmetric functions $f$. It is
not difficult to prove the analytic continuation of $\mathcal G_\lam$
for root systems of type A  using the Cayley-Capelli
type identities with the Dunkl operators
\cite{gz-br2}. (This will then give an refinement of
Theorem \ref{an-cont} showing that
$\Gamma_a(\frac {\delta-1 +\iota+2b}2 +\frac a2(r-1)+1)|\SH(t)|^{\delta}$
has an analytic continuation.)
 However it seems that to prove 
the result on Dirac distribution
some more works on the Laplace transform
\cite{Baker-Forrester-Duke}
on root system  of type A have to be done.
\end{rema+}

\def\cprime{$'$} \newcommand{\noopsort}[1]{} \newcommand{\printfirst}[2]{#1}
  \newcommand{\singleletter}[1]{#1} \newcommand{\switchargs}[2]{#2#1}
  \def\cprime{$'$} \def\cprime{$'$}
\providecommand{\bysame}{\leavevmode\hbox to3em{\hrulefill}\thinspace}
\providecommand{\MR}{\relax\ifhmode\unskip\space\fi MR }
\providecommand{\MRhref}[2]{%
  \href{http://www.ams.org/mathscinet-getitem?mr=#1}{#2}
}
\providecommand{\href}[2]{#2}

\end{document}